\documentclass[12pt,letterpaper]{article}
\usepackage[utf8x]{inputenc}
\usepackage{float}
\usepackage{quoting}
\usepackage{mdframed}
\usepackage{dsfont}
\usepackage{amsthm}
\usepackage{thm-restate}

\usepackage[
letterpaper,
top=1in,
bottom=1.2in,
left=0.7in,
right=0.7in]{geometry}

\usepackage{bm}
\usepackage{hyperref}
\hypersetup{colorlinks=true, urlcolor=blue, linkcolor=blue, citecolor=green, }
\usepackage{prettyref}

\usepackage{amssymb}
\usepackage{amsfonts}
\usepackage{latexsym}
\usepackage{amsmath}
\usepackage{cleveref}

\usepackage[normalem]{ulem}
\usepackage{todonotes}
\hypersetup{colorlinks=true, urlcolor=blue, linkcolor=blue, citecolor=magenta}

\newtheorem{thm}{Theorem}

\renewcommand{\phi}{\varphi}

\newcommand{\N}{\mathbb{N}}

\newcommand{\R}{{\mathbb R}}
\newcommand{\C}{\mathbb{C}}

\newcommand{\ket}[1]{{|#1\rangle}}
\newcommand{\bra}[1]{{\langle#1|}}

\begin{document}

\title{Knot theory and quantum computing}
\author{Robin Gaudreau and David Ledvinka}
\date{\today}

\maketitle 
\abstract{This paper explores the interactions between knot theory and quantum computing. On one side, knot theory has been used to create models of quantum computing, and on the other, it is a source of computational problems. Knot theory is often used to introduce topological idea to people without a formal mathematical background, and we are building on this tradition to discuss some of the deeper ideas of quantum computing.}

%Reference for intro to theoretical computer science

\section*{Introduction}

\paragraph{Disclaimer.} This text is aimed at non-physicists and does not concern itself with the feasibility of building a quantum computer. It treats quantum computation as a purely mathematical model. Conversely, much of the subtleties of knot theory have been omitted as an attempt to streamline the text. 

\subsection*{Knot theory}
Knot theory is the mathematical study of an idealized model of knots. It primarily uses algebraic and geometric techniques to study topological objects.  A \emph{knot} is a smooth embedding $\kappa: S^1\hookrightarrow \R^3$, considered up to isotopy (smooth continuous deformations).

\begin{figure}[htpb]
    \centering
    \includegraphics[width=0.2\textwidth]{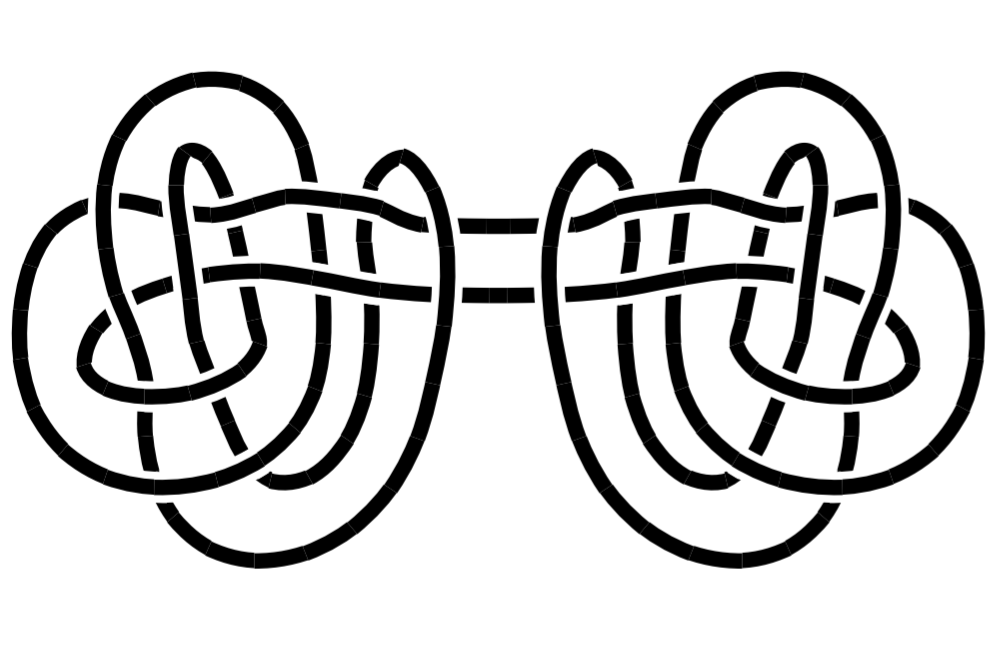}  \hspace{0.02\textwidth} 
    \includegraphics[width=0.2\textwidth]{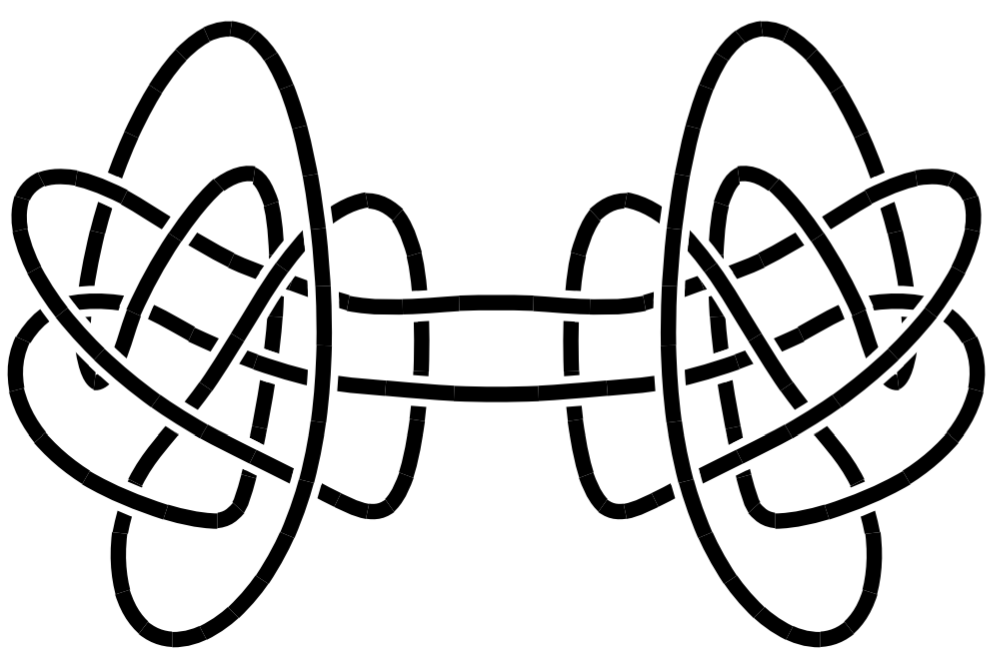} \hspace{0.02\textwidth} 
    \includegraphics[width=0.2\textwidth]{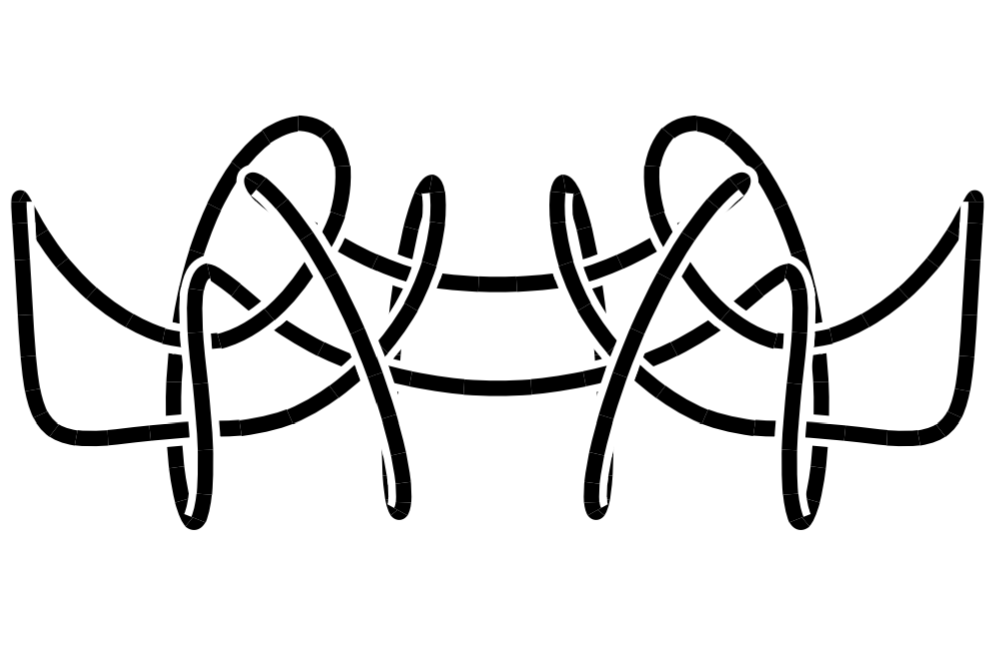} \\
    \includegraphics[width=0.2\textwidth]{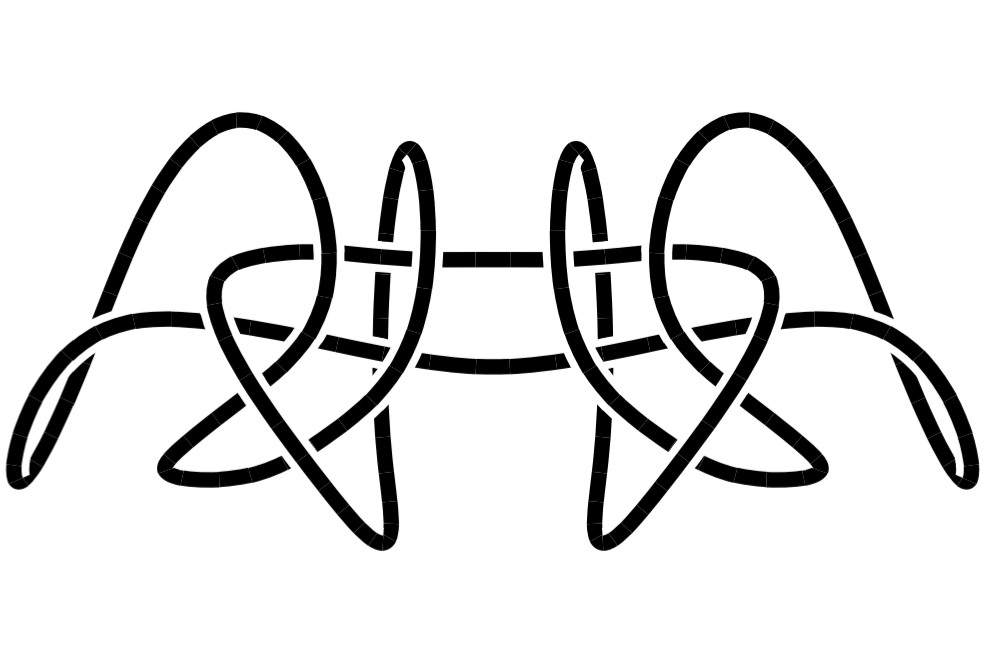} \hspace{0.02\textwidth}  
    \includegraphics[width=0.2\textwidth]{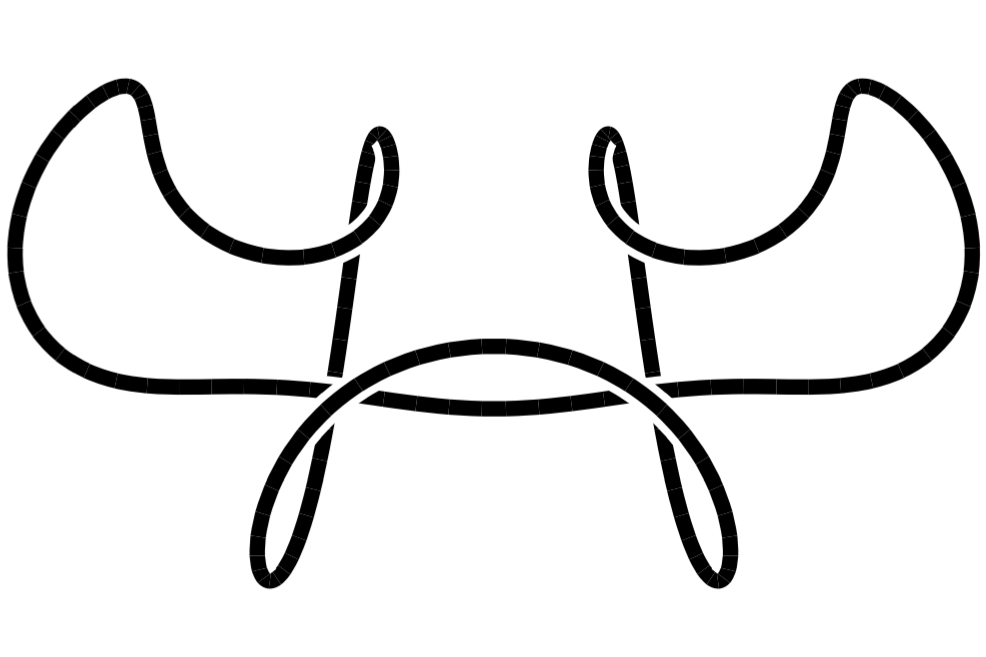} \hspace{0.02\textwidth} 
    \includegraphics[width=0.2\textwidth]{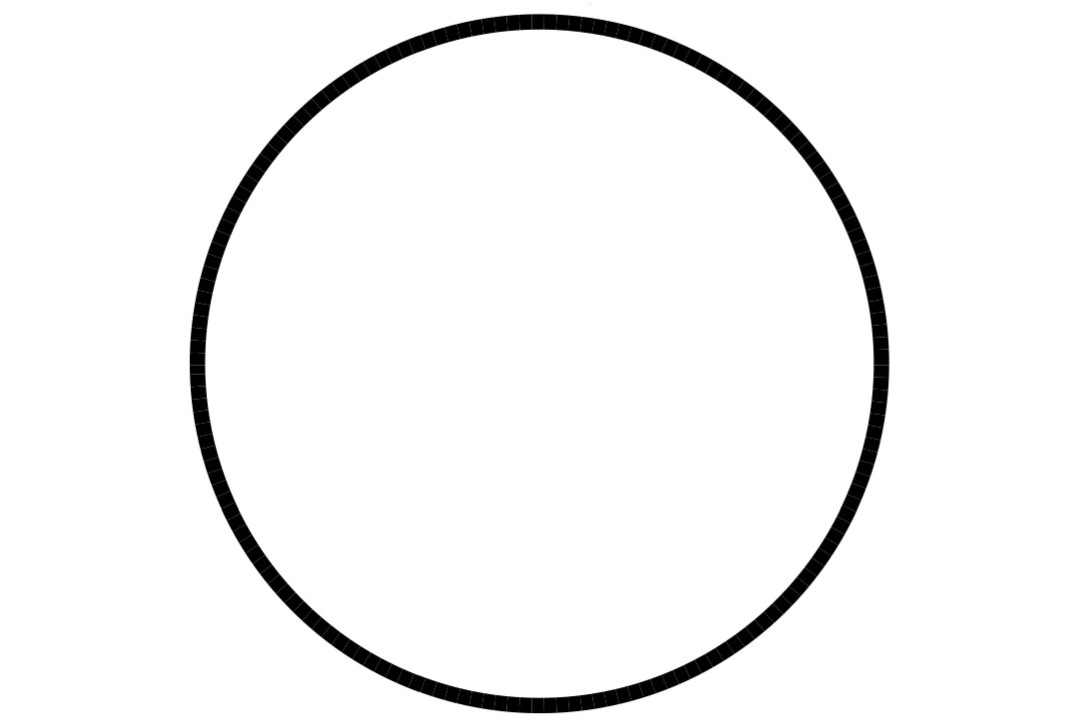} 
\caption{Various diagrams of the unknot, from \cite{knotplot}.}
\label{diagrams}
\end{figure}

Informally, one thinks of a knot as any of the shapes that can be made by a perfectly elastic string which has been tangled in space and whose ends are glued together. We consider two such shapes to be the same knot if we can we can manipulate one in space to get the other without pulling the string through itself. They are represented by  \emph{knot diagrams}, which are curves in the plane with transverse self-intersections called crossings, which are represented as in Figure \ref{diagrams}. Knot diagrams which are related by a finite sequence of the operations in Figure \ref{rms} represent isotopic circle embeddings. From an embedding of a circle, one creates knot diagrams by projecting it to a plane, and from a knot diagram, one recovers a smooth embedding of a circle by resolving the crossings in the direction perpendicular to the plane in which the curve is drawn. Therefore, knot diagrams encode all the information needed to recover an isotopy class. 

One of the oldest problems in knot theory is the question of obtaining an exhaustive and non-redundant list of knots. As there are countably infinitely many knots, these lists, called knot tables, are usually listed by the minimal number of crossings in a planar projection of the knot. The first such enumeration was Taits' ``First Seven Orders of Knottiness", which was inspired by Lord Kelvin's theory that atoms were small knotted vertices in aether, and that their properties came from the topology of the knots. While that model was both quite inaccurate, and gathered little attention, there is poetic justice in the way this idea motivated tools which are now used in the study of particle physics. For more historical details, the reader can consult \cite{Ga2016}, and more notably the references therein.

\subsection*{Quantum computing}
Classical computing can be summarized as encoding information in binary, modifying it using (usually deterministic) rules, and outputting a binary answer.

In contrast, quantum computing stores information as tensor products of elements of $\C \operatorname{P}^1=\{[z_0:z_1)]: z_i\in \C, \text{ and } z_0 \not=0 \text{ or } z_1 \not = 0, \ [z_0:z_1]\sim [\lambda z_0: \lambda z_1] \ \forall \ \lambda \in \C^*\}$. The smallest unit of information for a quantum computer is the \emph{qubit}. It is assigned an element of $\C \operatorname{P}^1$ called its \emph{state}, which is traditionally written using Dirac's bra-ket notation for ease of manipulation by \emph{gates} encoded as projective unitary matrices. Finally, the output is given by projecting to an orthonormal basis in a process called \emph{measurement}. 

\paragraph{Example.} A qubit in the state $[\alpha:\beta]$ is written $\alpha \ket 0+ \beta \ket 1$. It measure to $\ket i$ with probability $|(\alpha \bra 0 + \beta \bra 1 ) \ket i|^2$. An example of a unitary transformation is the Hadamard gate, $H=\frac{1}{\sqrt 2}\left[\begin{matrix}1 & 1 \\ 1 & 1\end{matrix}\right] $. It acts on a single qubit, such as $H \ket 0 = \frac{1}{\sqrt 2} (\ket 0 + \ket 1)$.

Theoretically, a qubit is just a point on the projective line. Practically, a qubit is a transistor which can encode a state $\alpha \ket 0 + \beta \ket 1$ where $|\alpha|^2 + |\beta|^2 = 1$. A quantum circuit is then a sequence of physical forces acting on a collection of qubits. Qubits whose states are non-trivial tensor products are said to be \emph{entangled}. They are obtained by multiplying two unentangled qubits by some $4\times 4$ matrices.

\subsection*{Structure of the paper} In Section \ref{braid}, the algebraic background needed to understand the applications of knot theory in regards to quantum computing is presented. This consists of an introduction to knot invariants from the study of unitary representations of braid groups. Then, Section \ref{computer} explains the basics of a model of perturbation-resistant quantum computer and surveys the work that has been done to emulate how such a computer would implement a quantum circuit. In Section \ref{complex}, the converse interaction between knot theory and quantum computing is explored by comparing the classical and the quantum computational complexity of some knot invariants. Section \ref{conclusion} concludes with some open problems. 

\section{Braid group representation} \label{braid} 
\begin{quotation} \raggedleft
  \emph{God created the knots, \\ all else in topology is the work of mortals.}\\
    -- Dror Bar-Natan, modified from Leopold Kronecker.
\end{quotation}

A powerful source of information about the isotopy class of a knot diagrams is representation theory: the study of ways to associate a matrix to each element of a group. This gives the first interaction between knots and quantum computing, as knots can be associated to unitary matrices, although in a highly non-unique way.

\subsection{Braid notation}

A formal and authoritative source of information on knot theory is the frequently revised book by Burde and Ziechang, \cite{BZ2013}. The following draws from that text.

\subsubsection{Alexander's theorem} \label{sect - Alexander's theorem}

Braids are a convenient notation for knots and links. For our purpose, consider a \emph{braid diagram} on $n$ strands to be a collection of $n$ curves in the plane, oriented monotonically in the $x$-direction and crossing each other transversely. A braid is the equivalence class generated by a braid diagram up to the bottom two moves in Figure \ref{rms}.

\begin{thm}[Alexander, 1923]
Any knot can be represented as the closure of a braid.
\end{thm}

\begin{figure}[ht]
    \centering
    \includegraphics[width=0.7\textwidth]{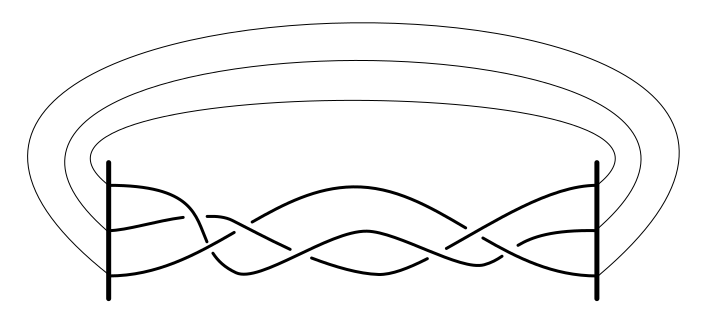} 
   \caption{The braid $\sigma_1 \sigma_2^{-1}\sigma_1 \sigma_2^{-1}\sigma_2\sigma_1^{-1}\sigma_2\in B_3$ closing to a link.}
   \label{closure}
\end{figure}

Here, the closure is the operation that is realised by gluing the thin lines in Figure \ref{closure} to the bold braid diagram. The knot or link obtained by closing $\beta\in B_n$ in such a way is called $\hat \beta$. By Vogel's algorithm \cite{Vo1990}, a knot diagram can be put into a braid-like diagram in quadratic time.

\subsubsection{Braid groups}

Fix $n\in \N$, $n\ge3$. The braid group on $n$-strands, denoted $B_n$ is defined as the set of finite words in the alphabet $\{\sigma_1^{\pm1}, \ldots \sigma_{n-1}^{\pm 1}\}$, subject to the relations: 

\begin{enumerate}
      \item Invertibility, $\sigma_i\sigma_i^{-1}= \sigma_i^{-1}\sigma_i=1$, where $1$ is also written as the empty word.
    \item Far commutativity, $\sigma_i^\epsilon\sigma_j^\eta=\sigma_j^\eta\sigma_i^\epsilon$ for all $|i-j|>2$, and $\eta, \epsilon \in \{-1, +1\}$;
    \item Yang-Baxter relation, $\sigma_i\sigma_{i+1}\sigma_i= \sigma_{i+1}\sigma_i\sigma_{i+1}$.

\end{enumerate}

The first relation justifies the use of the name `group' for this object. The generator $\sigma_i$ correspond to the $i$th strand from the top crossing on top of the $i+1$st strand. Its inverse is the move where the $i$th strand crosses behind the next. A word in that alphabet is read from left to right, follows the orders of the crossings in a braid diagram also from the left to the right. 

Making the generators their own inverses yields the map $ B_n\to S_n$, which maps the braid group on $n$ strands to the symmetric group on $n$-elements, also known as the group of permutations. A braid $\beta \in B_n$ such that $\operatorname{ab}(\beta)=e$ is a \emph{pure braid}. By this definition as the kernel of a group homomorphism, the set of pure braids, $P_n$, also forms a group. 

\subsubsection{Markov's theorem} %check the Birman-Menasco no-stabilization Markov theorem

Braids are helpful for enumerating knots because the braid group admits orderings. However, even braid words which do not represent the same group element can have isomorphic closure. 

\begin{thm}[Markov]
Let $\alpha\in B_n$ and $\beta\in B_m$, $n\le m$ be braids such that $\hat \alpha =\hat \beta$. Then, $\alpha$ and $\beta$ are related by a finite sequence of the following moves: 
\begin{enumerate}
    \item Right stabilization, $\alpha\mapsto \alpha \sigma_n\in B_{n+1}$, and
    \item Conjugation, $\beta \mapsto \gamma \beta \gamma^{-1}$ for some $\gamma\in B_m$. 
\end{enumerate}
\end{thm}

This explains why checking whether two braids represent the same knot is heuristically harder than a word isomorphism problem in a finitely generated group. It is well-known that a general word isomorphism problem in a fixed group is not solvable by a deterministic algorithm in polynomial time. Because of the stabilization operation, braids can only be compared as elements of a formal closure of the braid groups, $B_*=\cup_{n=2}^\infty B_n$.

\subsection{Invariants}
A function on knot diagrams which assigns the same value to all representatives of a knot is called a \emph{knot invariant}. Knot theorists say that some invariant $f$ dominates $g$ if there exists a pair of knots $K_1$ and $K_2$ such that $g(K_1)=g(K_2)$, but $f(K_1)\not=f(K_2)$. This inequality is proof that $K_1$ and $K_2$ are not isomorphic, and this is how knot invariants are used. The strongest possible knot invariant is one which takes a different value for each knot. It is rewarded with the title of \emph{classifying invariant}, but from a computational point of view, either distinguishing or computing the values of such an invariant has to be at least as hard as distinguishing knots themselves. Therefore, one has to settle for knot invariants which are computable and which take values that are easily distinguishable, such as complex numbers or Laurent polynomials.

\subsubsection{Skein formulas}

Some knot invariants are computed from an embedding in $\R^3$. Here, let's focus on computing from a knot diagram. Given as input one representative of a knot or link, it is possible to check that the function is invariant over the isotopy class by applying Reidemeister's theorem.

\begin{figure}[htbp]
    \centering
    \includegraphics[width=0.5\textwidth]{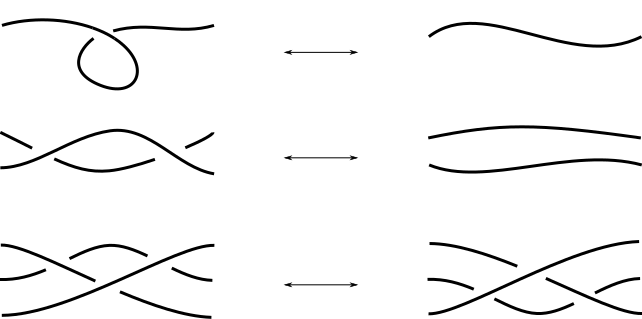}
    \caption{The three types of Reidemeister moves.}
    \label{rms}
\end{figure}

\begin{thm}[Reidemeister]
Two oriented link diagrams represent isomorphic links if and only if they are related by a finite sequence of the moves depicted in Figure \ref{rms}.
\end{thm}

\begin{figure}[ht]
    \centering
    \includegraphics[width=0.75\textwidth]{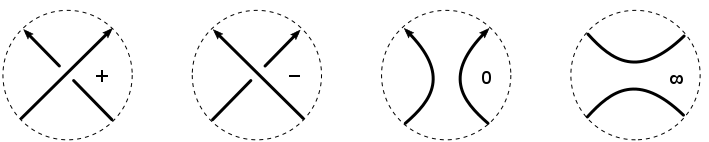} 
   \caption{The different ways to connect four nearby points in a knot diagram.}  \label{skein}
\end{figure}

To construct a knot invariant, look at the diagram as the sum of its parts. That is, define a function $f$ on the set of knot diagrams to be such that 

$$ Af(K_+)+ Bf(K_-)=Cf(K_0)+Df(K_\infty),$$

where the diagrams $K_+, K_-, K_0$, and $K_\infty$ are diagrams which agree everywhere except for the neighbourhood of one crossing in which they look like the sub-diagrams of Figure \ref{skein}. Then, choose the values of $A, B, C, D$ such that the function $f$ is invariant under those moves. \label{polyinv}

One such solution is the homflypt polynomial, which is given by the following skein relation:
$$ aP(K_+)|_{(a,z)}-a\-P(K_-)|_{(a,z)}=zP(K_0)|_{(a,z)},$$
with the condition that if $U$ denotes the unknot (the knot which admits a planar diagram with no crossings), then $P(U)=1$. Except for some special values of $(a,z)$, computing this invariant requires considering $2^c$ diagrams when starting with a knot that has $c$ crossings. This is an exponential time algorithm, far slower than the theoretical computer science ideal of \emph{efficient}, a label only applied to polynomial time algorithms. A quicker way to obtain information about the homflypt polynomial of a knot is explained in Section \ref{vertigan}. 

\subsubsection{Coloured Jones polynomials}

The Jones polynomial is a 1-variable specialisations of the homflypt polynomial. It has other generalizations called the coloured Jones polynomials which were originally defined by computing the invariant of a link formed by replacing the knotted curves with $n$ parallel copies of them, yielding the \emph{cable} of the link. Recently, \cite{HL2018} gave an extensive analysis of an algorithm which computes the coloured Jones polynomials from the original knot. They found that the computation time grew both with the number of crossings and the number of colours.

\subsubsection{Gassner matrices} %make sure this is the right invariant for this setting.

A group representation is a homomorphism $\Phi$ from a group $G$ to the automorphism group of a vector space, $\operatorname{Aut}(V)$. This means in particular that for any $g, h\in G$, $v\in V$, $\Phi(gh)(v)=\Phi(g)\Phi(h)(v)$. Given a basis, automorphisms of finite vector spaces are usually written as matrices. A representation is said to be \emph{unitary} if it maps every element of the group to a unitary matrix, that is, one whose inverse is its own conjugate transpose. The usual Jones polynomial is computable as a trace of a unitary representation of the pure braid group. More details about this constructions are found in \cite{BN2014}. For technical reasons, it can only be called a representation for pure braids, but the corrections for the general braid group are computationally inexpensive and theoretically un-illuminating. This representation is defined first on generators by mapping $\sigma_i \in P_n$ by
$$\sigma_i\longmapsto I_{i-1} \oplus\left(\begin{matrix}1-t & t \\ 1 & 0\end{matrix}\right)\oplus I_{n-i-1}\in \operatorname{Aut} (\C^n).$$ 

Then we can extend this map to the entire braid group by sending any word $\sigma_{i_1}\sigma_{i_2} \dots \sigma_{i_k} \in B_n$ to the product of the matrices assigned to its generators. Then to compute the Jones polynomial of any knot, we can first represent it as a braid as described in section \ref{sect - Alexander's theorem}, compute the representation of this braid, and then compute a certain trace of this matrix. In \cite{AJL2006}, it is shown that such a computation is efficiently implementable as a quantum circuit.

\section{Topological quantum computers} \label{computer}

\begin{quotation} \raggedleft
\emph{The final test of every new mathematical theory is its success\\ in answering pre-existent questions that the theory was not designed to answer.\\}
-- David Hilbert, 1926.
\end{quotation} 

Topological quantum computation is a fault-tolerant model of quantum computation which seeks to encode the state of a computation in topological data of a system as opposed to local data, which makes it naturally resistant to error caused by perturbation. It uses theoretical 2D quasi-particles called non-abelian anyons (which are believed to exist \cite{FKLW2003}) which have the property that exchanging them in space can cause arbitrary unitary transformations to be applied to their states, in contrast with fundamental 3D particles such as fermions or bosons which only experience a phase shift of -1 or 1 respectively when exchanged. Interestingly this model is connected to braids and the Jones polynomial.  In this section we describe a version of topological quantum computation which allows one to implement algorithms designed for the traditional quantum circuit model. We follow the approaches found in \cite{FS2018} and \cite{FKLW2003}, and more details can before found in each of these texts.

\subsection{Initialization}

A topological quantum computation is initialised by creating pairs of Fibonacci anyons (a certain kind of non-abelian anyon) from the vacuum. When two Fibonacci anyons -- from now, on simply anyons -- are fused they have the chance to either annihilate to the vacuum or combine into one new anyon. However the net charge of a system of anyons is always conserved, which has the consequence that a pair of anyons created from the vacuum will always annihilate when fused. During the computation we will exchange the position of these particles in space which will cause their states to evolve by a unitary operation and hence other fusion outcomes become possible. To mimic the quantum circuit model we will further group the pairs of anyons into groups of four (containing two pairs) which will represent one qubit. At the end of our computation we will fuse a fixed pair of anyons in each qubit to determine the result of our computation. In the event that they annihilate we consider the qubit to being measured as a $\ket{0}$, and in the event that they combine we consider it to be a $\ket{1}$. Since initially our anyon pairs will always annihilate, our qubits always start in the  $\ket{0}$ state.

\subsection{Computation}

\begin{figure}[ht] 
    \centering
    \includegraphics[width=0.5\textwidth]{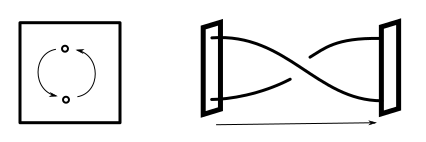}
    \caption{Swapping two particles in the plane and the corresponding braid in 2+1 spacetime.}
    \label{anyons}
\end{figure}

We now perform a computation by exchanging the positions of the anyons in space. When we do this a unitary operator is applied to their state which can change both the possible outcomes of fusion and the probabilities that each fusion outcome occurs. Swapping adjacent anyons in the plane braids the path of the anyons in 2+1 spacetime as seen in figure \ref{anyons}. Hence performing a sequence of swappings among all the anyons in the system traces out a braid in 2+1 spacetime, where each crossing corresponds to an exchange of anyons and hence an application of a unitary operator. Braids which entangle anyons from two different qubits literally produce entangled quantum states, while braids which act upon anyons within the same qubit correspond to one-qubit operators. It can then be shown \cite{FLW2002} that we can design a braid which approximates any unitary operator to arbitrary precision, in particular the common quantum gates used to form a universal gate set in the circuit model. Thus we can implement any quantum circuit by composing the braids corresponding to the gates used in the circuit in the order that they are applied. 

Since a topological quantum computation is solely dependent on the topology of the braid, the computation will be unaffected by local perturbations of the system, which is in contrast to other proposed implementations of quantum computation where data is stored locally, and hence is very sensitive to local changes. 

\subsection{Measurement}
Once the computation is performed we measure the states of the anyons by fusing them together to determine the result of the computation. Fibonacci anyons are called  anyons because the amount of possible fusion outcomes in a system of $n$ anyons is given by the $n$th Fibonacci number. But when we are implementing a quantum circuit only some of these outcomes will correspond to basis states in the quantum circuit model, which we will call computational states. Hence we will want to be careful to design braids which only lead to fusion outcomes that correspond to one of these states (alternatively we could allow for error). When we have this guarantee we can measure the value of a qubit simply by fusing a fixed pair of anyons within the qubit. If the pair of anyons annihilate then we consider the qubit to have a value of $\ket{0}$ and if the pair of anyons combine we consider the qubit to have a value of $\ket{1}$. After measuring each of the qubits we will have a string of bits which we will take to be the result of the computation.

\begin{figure}[ht]
    \centering
    \includegraphics[width=0.6\textwidth]{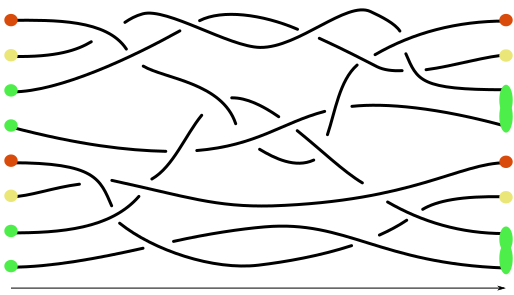}
    \caption{An example of a quantum circuit on 2 qubits acting by braiding together two quartets of anyons. The merging on the right hand side denotes fusion.}
\end{figure}

When the braids are designed properly, every individual qubit will have always have a net charge of $0$ meaning that it will annihilate if all four anyons in the qubit are fused. But pairs of anyons do not have to maintain this property so that fusing two anyons can result in them annihilating or combining. This is why only one pair of anyons has to be fused in order to determine the fusion outcome in each qubit, because the fusion outcome of the other pair will be entirely determined by the first. If the first pair annihilated then the second pair must annihilate. If the first pair combined into an anyon, then the second pair must do so as well, and when these two new anyons are fused they will annihilate. But because of this dependence, braiding the anyons from the second pair in a qubit will still effect the fusion outcomes of the system despite not being directly measured.

\subsection{Relationship with the Jones Polynomial}
Given a braid implementing a certain topological quantum computation, it turns out that a formula for computing the probability that the computation returns the all $0$ string is equivalent to an evaluation of the Jones polynomial of knot which is a certain closure of this braid. This is described in detail in \cite{FS2018}. This has implications for quantum complexity theory and the complexity theory of knot invariants. Since we believe calculating the exact probability that a quantum circuit outputs a certain value is a hard problem, this is further evidence that exact evaluations of the Jones polynomial is very hard. Since we believe even approximating this probability should be difficult for classical computers, this suggests that approximating evaluations of the Jones polynomial is a candidate for a problem which can only be solved efficiently by quantum computers. %(In BQP but not in P)

\iffalse
\subsection{Universal gate set}

Realisations of quantum gates for this model of topological quantum computing have been found using exhaustive search algorithms. A mathematician may want to think of this as a physical instance of an automorphism group, while a physicist would enjoy looking at this as a mathematical model for the physical system. \note{References from each approach}

\subsubsection{Another model of topological quantum computer}

To give an explicit result which uses the representation presented above, consider the result from \cite{FLW2002} that there exists a representation of $B_{2n}$ which is dense in $\operatorname{PSU}(H)$ for $n\ge 3$, where $H$ is a 2$n$-dimensional Hilbert space which can be interpreted as the register of a quantum computer. The restriction to the values of $r$ is due to the other numbers corresponding to roots of unity which do not generate a dense set in $\C$. The group $\operatorname{PSU}(H)$ is the set of unitary matrices with determinant $\pm 1$, modulo multiplication by scalar multiples of the identity matrix. 

\subsubsection{The one-clean qubit model}
\cite{JW2009} \problem{This was mentioned above.}
\fi
\section{Computational complexity of the homflypt polynomial} \label{complex}

\begin{quotation} \raggedleft
   \emph{ Using this exponentially large computational space, it is possible, at least in principle, for quantum computers to efficiently solve classically difficult problems. }\\    
-- Bernard Field and Tapio Simula, 2018.
\end{quotation}
%\paragraph{Homflypt determines the Jones and Alexander.} 
Recall the homflypt polynomial-valued invariant defined in Section \ref{polyinv} and that the (2-coloured) Jones polynomial of a knot $K$, is determined by $P(K)$. It is shown in \cite{JVW1990} that the Jones polynomial is essentially as difficult to compute as the homflypt polynomial. Throughout this section, the computation time is considered to be a function of the crossing number. 

\subsection{Classical computation}
There is a subtle difference between the existence of a mathematical formula for a quantity and the existence of an algorithm to compute it. Without going into the subtleties of works mentioned, this section reviews some of the work that has been done about the computational complexity of some polynomial knot invariants.

\subsubsection{Vertigan's algorithm} \label{vertigan}

As computing the entire homflypt polynomial is slow, we can compute individual coefficients efficiently and classically. Such invariants are said to be of \emph{finite type} since they can be computed on knots with some finite amount of information missing, up to a point at which they vanish. 

In this case, write the homflypt polynomial as a polynomial in $z$ with coefficients depending on $a$. If $P(K)|_z= \sum_i c_i(a) z^i$. Then, there is an algorithm to compute $c_i$, that is presented in \cite{Pr2017} and attributed to Vertigan that yields the following result.

\begin{thm}
The (classical) time complexity of computing the coefficient $z^{2i+1-k}$ of the homflypt polynomial of a $k$-component link is bounded by a polynomial in the number of crossings of the link of degree linear in $i$.
\end{thm} 

Recursive formulas for the coefficients allow us to consider a table of only $n^2$ knots with decreasing complexity, but we may need to use the homflypt polynomial of these knots exponentially many times. Saving those values to be re-used requires a polynomial amount of memory. See \cite{Le2018} for a Mathematica program which is based on this algorithm.

\subsubsection{Classical computational complexity}
Since the definition of the homflypt polynomial fits on a single line, it could be expected that there is a fast algorithm to compute it that nobody has just been clever enough to find. However, its computational complexity is bounded below by that of the Jones polynomial.

\begin{thm} [Jaeger, Vertigan, Welsh, 1990]
Determining the Jones polynomial of an alternating link is \#P-hard.
\end{thm}

This result is significantly more general than it sounds, because not only roughly half of the links admit an alternating diagram (that is, a diagram in which travelling along any component, one crosses crossings in an alternating over-then-under fashion), but the computation of invariants using skein formulas is simpler on alternating diagrams.  

\begin{thm}[Kuperberg, 2009]
Fix a link diagram $L$, a principal root of unity $t=e^{2\pi i/k}$, $k>6$ or $k=5$, and $0<a<b$. Given that $|V(L)|_t|<a$ or $|V(L)|_t|>b$, it is \#P-hard to decide which inequality holds.
\end{thm}

\#P is the class of enumeration problems associated with NP. This problem is reduced to an enumeration problem is by an algorithm which considers only using paths in the skein tree which contribute non-trivially to the final polynomial. 

%The natural relationship between the Jones polynomial and matrix multiplication hints that polynomial quantum time corresponds to exponential classical time.

\subsection{Quantum computation}

The quantum algorithms that compute the coloured Jones polynomials do so only at specific values since quantum circuits can only encode matrices with entries in $\C$, and not with entries in the more general ring $\mathbb Z [q^{-1/2}, q^{1/2}]$. To get an idea of how far quantum computing does have a chance to go, one should think of the quantum circuits described as sequences of unitary transformations as comparable to describing classical computing as a Turing machine acting on a strip of symbols. It is well knows that this is equivalent to algorithms in natural languages, but creating a machine that works exactly like Turing's theoretical model is tedious.  More discussion of the limitations of the quantum circuit model can be found in \cite{BBBV1997}. This yields hope that there is much room for improvement in this field of research, but unfortunately muddies attempts to compare the algorithms written for a computer algebra system with those written for a quantum computer. Moreover, a quantum computer cannot be asked to output the exact numerical answer with probability 1. The best it can do is provide an approximate answer, or in the case of decision problem,  yield a correct answer with high probability.

\subsubsection{Additive approximation}

Given a function $f: \C\to [0,1]$, a quantum algorithm $A$ gives an \emph{additive approximation} of $f$ if there exists $\epsilon >0$ such that for any $z\in \C$, $|f(z)-A(z)|<\epsilon$ with high probability. 

\begin{thm}[Aharonov, Jones, Landau, 2006] There exists a quantum algorithm which, for a braid $\beta\in n$ of length $m$, gives an additive approximation to $V(\hat \beta)|_{e^{2\pi i/k}}$.
\end{thm}

The additive approximation can be used to distinguish knots if the values can be bounded away from each other. The computation can be halted as soon as the polynomials are shown to differ at at least one point. This is something which may be useful when dealing with very large knots, where a precise analysis is beyond the memory resources of current computers. 

For example, if two quantum circuits $C_1$ and $C_2$ are claimed to output the same measurement, the likeliness of this statement being true can be tested by evaluating the Jones polynomial of the braids obtained by representing the two circuits as the braids, respectively $b_1$ and $b_2$, that would be created by the anyons in a topological computer running those circuits. If the braids are isomorphic, then the outputs have to agree, and so do the Jones polynomials $V(\hat b_1)$ and $V(\hat b_2)$. 

\subsubsection{Quantum complexity of computing the Jones polynomial}

According to \cite{Wa}, the possibility that NP is contained in BQP is a motivation for research in quantum computing. The complexity class {BQP} (bounded error quantum polynomial time) contains problems that can be solved efficiently on a quantum computer in polynomial time, and correctly with probability $2/3$.  

In evaluating the polynomial a quantum computer creates uncertainty and needs to be run multiple times. Since the Jones polynomial is a meromorphic function, the whole polynomial is determined by its value at a sufficient large number of primitive roots of unity. Moreover, knowing that its coefficients are integers compensates for the error introduced by the approximation. 

It turns out that estimating even the regular Jones polynomial is a rather general problem for quantum computing. A common theoretical construct used to analyse the relationship between complexity classes is the \emph{oracle}. It is a single gate which computes a function without the computer or the mathematician needing to know how this function is defined. 

\begin{thm}[Bordewich, Freedman, Lovasz, Welsh, 2005]
Let $A(\beta,z)$ be an oracle that computed the Jones polynomial of a braid $\beta$ at $z$, then \emph{BQP} $\subset \operatorname{P}^A$.
\end{thm}

In other words, giving an additive approximation of the Jones polynomial is a BQP-complete problem. For more information about the kind of problems that are in BQP, see \cite{Wa}.

\section{Conclusion} \label{conclusion}
%Approximations of the homflypt polynomial can be used to distinguish knots. Computing exponentially many of them (at different $k$-th roots of unity) can yield enough data to fit the polynomial (given that it has integer coefficients). So this does not speed up the computation of the whole polynomial, but in general, to tell two knots apart, the computation can be halted as soon as two values are determined to be bounded away from each other.

\begin{quotation} \raggedleft
   \emph{Slowly, the implications of the idea began to be understood. To begin with it had been too stark, too crazy, [...] then some phrases like `Interactive Subjectivity Frameworks' were invented, and everybody was able to relax and get on with it.} \\
   -- Douglas Adams, in \emph{Life, The Universe And Everything}, 1982.
\end{quotation}

Knot theory and quantum computing intersect in many ways that were beyond the scope of this introductory paper. It is only fair to a reader who made it this far to at least mention them. There is the result of \cite{WY2008} that a promise problem, choosing between an upper bound and a lower bound for the magnitude of the Jones polynomial, is QCMA-complete (the complexity class QCMA is denoted MQA in \cite{Wa}). The One Clean Qubit computation model, as explored in \cite{JW2009}, also employs braid representation, and finally, topological quantum field theory entangles quantum computing and knot invariants in a natural way. A concise reference for that last topic is \cite{FLW2002}. Let us now state, as promised, some open problems in quantum computing coming from knot theory.

\paragraph{1. Planarity of intersection sequences.} Given a knot diagram, assign to each crossing a number, and, travelling along the knot, write down the numbers as they are encountered and how (under of over) the crossing is crossed. The result is the intersection sequence, also known as Gauss code, of the knot. In general, words that look like intersection sequences do not come from knot diagrams. The planarity problem is to determine if a given sequence is the intersection sequence of a knot. There are many algorithms which yield a non-deterministic solution to the planarity problem. Is this problem in NP? What is its quantum computational complexity?

\paragraph{2. Quantum algorithms for finite type invariants.} Coefficients of the homflpt polynomial are not the only kind of finite type invariants. There appears to be a gap in the literature around quantum algorithms for other finite type knot invariants.

\paragraph{3. Geometry of quantum processors.} The model of topological quantum computer presented here, like many other quantum computers assumes that qubits are forced to live in processors consisting of a piece of plane. This reduces the number of elementary interactions between them. This can be changed by instead placing the qubits on a surface of a higher genus. The corresponding mathematical object is the fundamental group of the configuration space of $n$ points on a surface. Could such a computer be built? What kind of time savings can be done by augmenting the number of direct neighbours a qubit has?\\

\textit{Acknowledgements. } We would like to thank Dror Bar-Natan and Hans Boden for their continual support and the many energizing conversations. This work was financed by the National Science and Engineering Research Council of Canada, and completed as part of course work for \emph{Quantum Computing: Foundations to Frontier} given by Henry Yuen at the University of Toronto in the Fall of 2018.

\end{document}